# Fast Intrinsic Mode Decomposition and Filtering of Time Series Data

Louis Yu Lu

E-mail: louis-lu@rogers.com

Abstract: The intrinsic mode function (IMF) provides adaptive function bases for nonlinear and non-stationary time series data. A fast convergent iterative method is introduced in this paper to find the IMF components of the data, the method is faster and more predictable than the Empirical Mode Decomposition method devised by the author of Hilbert Huang Transform. The approach is to iteratively adjust the control points on the data function corresponding to the extrema of the refining IMF, the control points of the residue function are calculated as the median of the straight line segments passing through the data control points, the residue function is then constructed as the cubic spline function of the median points. The initial residue function is simply constructed as the straight line segments passing through the extrema of the first derivative of the data function. The refining IMF is the difference between the data function and the improved residue function. The IMF found reveals all the riding waves in the whole data set. A new data filtering method on frequency and amplitude of IMF is also presented with the similar approach of finding the residue on the part to be filtered out. The program to demonstrate the method is distributed under BSD open source license.

## 1 Introduction to Hilbert Huang Transform (HHT)

The real time series numeric data from natural phenomena, life science, social and economic systems are mostly nonlinear and non-stationary. The traditional analysis methods like Fourier transform and Wavelet transform presume the linear and stationary of the underlying system, the predefined function bases have little bearing on the physical meaning of the system, and are difficult to reveal the nature of the real data.

The recent researches tend to use adaptive function bases, Huang et al introduced the IMF as adaptive a posteriori function base in the form of Hilbert spectrum expansion, which has meaningful instantaneous frequency [1][2], The HHT is the Hilbert transform applied to the IMF components.

The HHT has broad application to signal processing and times series data analysis in the fields of health, environmental, financial and manufacturing industries [4].

### 1.1 Intrinsic Mode Function (IMF)

Dr Huang analyzes the requirement of meaningful instantaneous frequency on Hilbert transformation, and introduced the adaptive function base as IMF [1] [2], which satisfies the two conditions:

1) In the whole data set, the number of zero crossings must equal or differ at most by one;
2) At any point, the mean value of the envelope defined by the local maxima and the envelope defined by the local minima is zero.

The first condition is apparently necessary for oscillation data; the second condition requires the symmetric upper and lower envelopes of IMF, as the IMF component is decomposed from the original data; it is quite challenging to find the real envelopes because of nonlinear and non-stationary nature in the data. Only a few functions have the known envelopes, for example, the constant amplitude sinusoidal function.

All the IMF has the meaningful instantaneous frequency defined by the Hilbert transformation [3]. The different IMF components of a time series data reflects the oscillation of the data on different time scales, while the residue components reflect the data trend.

## 1.2    Empirical Mode decomposition (EMD)

To break up the original data into a series of IMF, Huang et al invented the EMD method [1] [4], the idea is to separate the data into a slow varying local mean part and fast varying symmetric oscillation part, the oscillation part becomes the IMF and local mean the residue, the residue serves as the input data again for further decomposition, the process repeats until no more oscillation can be separated from the residue of higher frequency mode.

On each step of the decomposition, since the upper and lower envelope of the IMF is unknown initially, a repetitive sifting process is applied to approximate the envelopes with cubic spline functions passing through the extrema of the IMF, the data function serves as the initial version of IMF, and the refining IMF is calculated as the difference between the previous IMF version and mean of the envelopes, the process repeats until the predefined stop condition satisfied. The residue is then the difference between the data and the improved IMF.

There are four problems with this decomposition method:

1)  The spline (cubic) connecting extrema is not the real envelope, some unwanted overshoot may introduced by the spline interpolation, the resulting IMF function does not strictly guarantee the symmetric envelopes [1]. Higher order spline does not in theory resolve the problem.

2)  The different stopping condition value results different set of IMF [1], picking up the right value is a more subjective matter, making the results unpredictable.

3)  The repetitive sifting process is time consuming, because the problem in item 2, more sifting does not guaranty better results. Some researches try to improve the performance, the study in reference [6] only speeds up the envelopes building process; problem 1 and 2 still exist.

4)  Some riding wave on steep edge of IMF are missing [1], the resulting IMF may not accurately reflect the true physical nature of the data.

## 2    Sawtooth transform decomposition method [5]

The author had introduced an efficient sawtooth decomposition method. The approach is to transforms the original data function into a piecewise linear sawtooth function (or triangle wave function), then directly constructs the upper envelope by connecting the maxima and construct lower envelope by connecting minima with straight line segments in the sawtooth space, the IMF is calculated as the difference between the sawtooth function and the mean of the upper and lower envelopes. The results found in the sawtooth space are reversely transformed into the original data space as the required IMF and envelopes mean. This decomposition method processes the data in one pass to obtain a unique IMF component.

The problem with the method is the leaking of high frequency ripples into the residue function. The following figures depict the first intrinsic mode decomposition result using the sawtooth method; the sample signal is the sum of two sinusoidal functions:

$$x(t) = 70\cos(\frac{\pi t}{150}) + 30\cos(\frac{\pi t}{15})$$

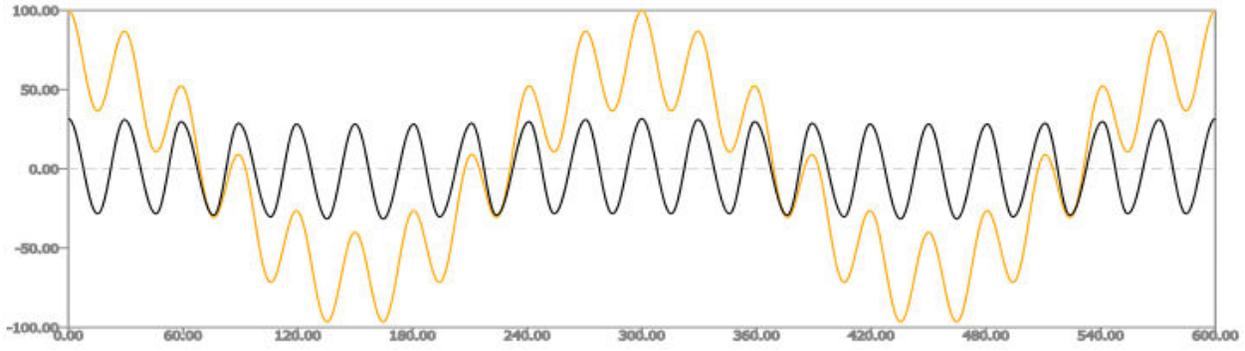

Figure 2.1 Original signals (in orange) and the first mode IMF (in black)

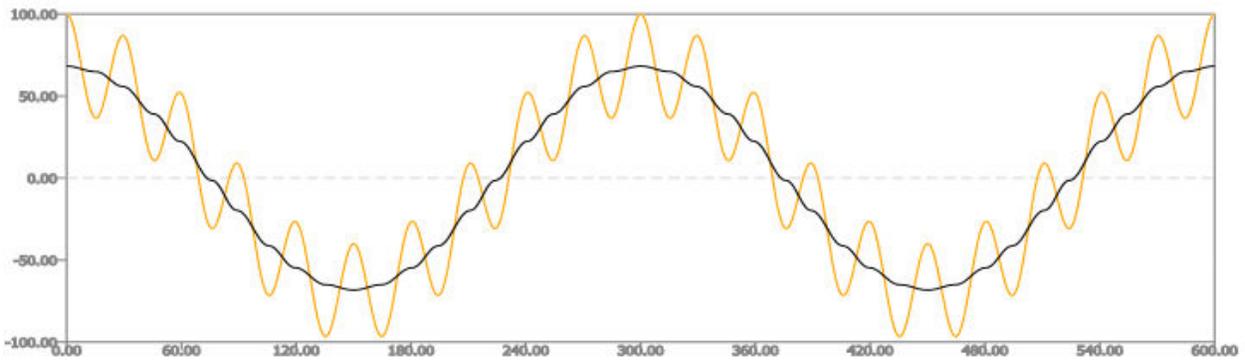

Figure 2.2 Original signals (in orange) and the first residue function (in black)

The figure 2.2 shows the undesired high frequency ripples leaked into the residue function, and they will end up into the low frequency IMF.

## 3   Fast intrinsic mode decomposition (Fast IMD)

The new method presented in this paper is to approximate the residue function through a fast convergent iterated process. The key is to easily find the initial residue function and quickly improve it through the iteration.

### 3.1 Initial the residue function and the IMF

Let *f(t)* be the original time series data and *f'(t)* the first derivative, *f'(t)* has *m* extrema:

$$R'(t_j) \qquad t_0 \le t_j \le t_{m-1}$$

The corresponding points on the data function with the same *t* values are inflection points:

$$R(t_j) \qquad t_0 \le t_j \le t_{m-1}$$

The initial residue function is built as the straight line segments passing through the inflection points and the initial IMF is calculated as the difference between the data function and the piecewise linear residue function. The initial IMF chosen this way can reveal the riding waves on the resulting IMF regardless of the steep edges of the data function.

An alternative way to choose the initial IMF is to use the data function **f(t)** directly, in this case the initial residue function is constant zero, with this initialization method, the resulting IMF will not reveal all the riding waves on the steep edges of the data function.

## 3.2 Improving the residue function and the IMF

With the current version of IMF, its **m-1** extrema found correspond **m-1** control points on the data function with the value $E(t_i)$; the coordinates of the control points are defined as:

$$P_i(t_i, E(t_i)) \qquad t_0 \leq t_i < t_{m-1}$$

Connecting these control points can construct a piecewise linear function; the key is to find the median points serving as the control points of the residue function.

The turning direction of two contiguous segments is defined as the vector product:

$$\alpha_i = P_{i-1}P_i \times P_iP_{i+1} \tag{3.1}$$

Since all the segments are in the same plane, the vector product has only non zero component on the perpendicular direction. The dot product used below can also be simply calculated as the product of values projected on the perpendicular direction.

For the current point $P_i$, if either of the previous or the following turning direction is different, the median point is calculated as the average of the control point and the point on the segment connecting the previous and following control points:

$$R(t_i) = \frac{1}{2}(E(t_i) + E(t_{i-1}) + (E(t_{i+1}) - E(t_{i-1}))\frac{t_i - t_{i-1}}{t_{i+1} - t_{i-1}})$$

$$\text{if } \alpha_{i-1} \bullet \alpha_i < 0 \text{ or } \alpha_i \bullet \alpha_{i+1} < 0 \tag{3.2}$$

Otherwise, both the previous and the following turning direction are on the same direction, the median point is identical to the control point:

$$R(t_i) = E(t_i) \text{ if } \alpha_{i-1} \bullet \alpha_i \geq 0 \text{ and } \alpha_i \bullet \alpha_{i+1} \geq 0 \tag{3.3}$$

The improved residue function is constructed as the cubic spline passing through the median points. The new version of IMF is calculated as the difference between the data function and the spline residue function. This method can also be used to find the median of any piecewise linear functions.

The process can be repeated by taking the new IMF as the current version to further improve the residue function and the IMF. The stop condition of the process can be the maximal iteration count, or the maximum point distance on successive residue functions falls in the given limit or the distance stopped decreasing.

## 3.3    IMF decomposition algorithm

The intrinsic mode components can be decomposed in the following steps:
1) Calculate the first derivative of the input data, and find the maxima and minima of the derivative. If the number of extrema is less than three, all the intrinsic modes are all found, stop the process;

2) Use the points on the data function corresponding to the extrema of the derivative to build straight line segments function as the initial residue function, calculate the initial IMF as the difference between the data and residue function;

3) Find the extrema on the current IMF and the corresponding control points on the data function, calculate the turning direction with formula (3.1), calculate the median points with formula (3.2) or (3.3), build the improved residue function with cubic spline passing through the median points and calculate the new IMF as the difference between the data function and the spline residue function;

4) If the maximal iteration count is reached, or the maximum point distance on successive residue functions is less than the given limit or the distance stopped decreasing, a component is successfully separated and proceed to the next step; otherwise go back to step 3);

5) Take the residue as input data again and go back to step 1).

### 3.4    Function extension and boundary points

To calculate the residue and IMF on the boundary points, the data series must be extended on both ends to ensure the continuity on starting and ending points. Based on the nature of the data, several possible extensions can be applied around the boundary points: even extension; odd extension or cyclic extension if the starting point and the ending point has the equal value.

The extension is only required for the control points, only two extra points are needed on each end.

1) Even extension: the starting point is considered as the first extrema $E(t_0)$, and the ending point as the last extrema $E(t_{m-1})$.

$$t_{-1} = t_0 - (t_1 - t_0) \qquad\qquad E(t_{-1}) = E(t_1)$$
$$t_{-2} = t_0 - (t_2 - t_0) \qquad\qquad E(t_{-2}) = E(t_2)$$
$$t_m = t_{m-1} + (t_{m-1} - t_{m-2}) \qquad E(t_m) = E(t_{m-2})$$
$$t_{m+1} = t_{m-1} + (t_{m-1} - t_{m-3}) \qquad E(t_{m+1}) = E(t_{m-3})$$

Even extension normally generates better result than other extension methods. For real data collected, better decomposition results can be obtained by starting and ending the sample data on the extrema points.

2) Odd extension: given the starting point $f(t_s)$, and the ending point $f(t_e)$.

$$t_{-1} = t_s - (t_0 - t_s) \qquad\qquad E(t_{-1}) = f(t_s) - (E(t_0) - f(t_s))$$
$$t_{-2} = t_s - (t_1 - t_s) \qquad\qquad E(t_{-2}) = f(t_s) - (E(t_1) - f(t_s))$$
$$t_m = t_e + (t_e - t_{m-1}) \qquad\qquad E(t_m) = f(t_e) - (E(t_{m-1}) - f(t_e))$$
$$t_{m+1} = t_e + (t_e - t_{m-2}) \qquad E(t_{m+1}) = f(t_e) - (E(t_{m-2}) - f(t_e))$$

3) Cyclic extension: the starting point and end point are considered as extrema and they have the same value $E(t_0) = E(t_{m-1})$.

$$t_{-1} = t_0 - (t_{m-1} - t_{m-2}) \qquad E(t_{-1}) = E(t_{m-2})$$
$$t_{-2} = t_0 - (t_{m-1} - t_{m-3}) \qquad E(t_{-2}) = E(t_{m-3})$$
$$t_m = t_{m-1} + (t_1 - t_0) \qquad\qquad E(t_m) = E(t_1)$$
$$t_{m+1} = t_{m-1} + (t_{m-1} - t_{m-3}) \qquad E(t_{m+1}) = E(t_2)$$

### 3.5 Comparison to other decomposition methods

Compare with EMD, the fast IMD method has the following advantages:

1) The method converges very quickly, satisfactory result can be obtained even with two iterations; for each step, only one spline is required for the residue function instead of two in EMD for upper and lower envelopes, it reduces the calculation time.

2) The stop condition is easy to choose. The shape of the resulted residue function and IMF is very stable, they stay the quite same shape even with excessive number of iterations, while the IMF mutates into equal amplitude frequency modulation function with EMD. Besides, the control points have doubled on the residue than on the envelopes, the effect of overshooting on the spline function is diminished.

3) The riding wave on the steep edge of the data is also revealed on the IMF, the resulted components reflect better the true nature of the data.

Compare with the sawtooth method, the fast IMD resolved the leaking problem of high frequency ripples into the residue function. Figures 3.1 and 3.2 show the new results relating to the flawed decomposition in figures 2.1 and 2.2.

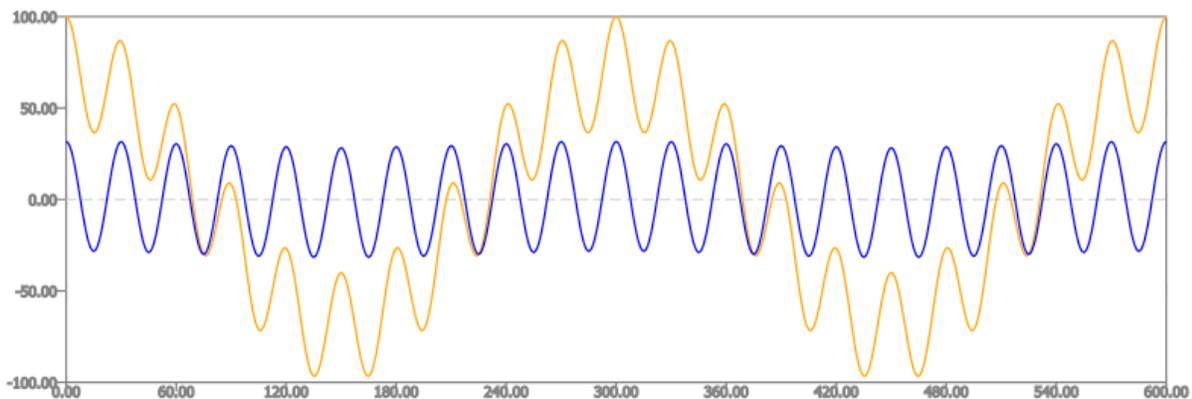

Figure 3.1 Original signals (in orange) and the first mode IMF (in blue)

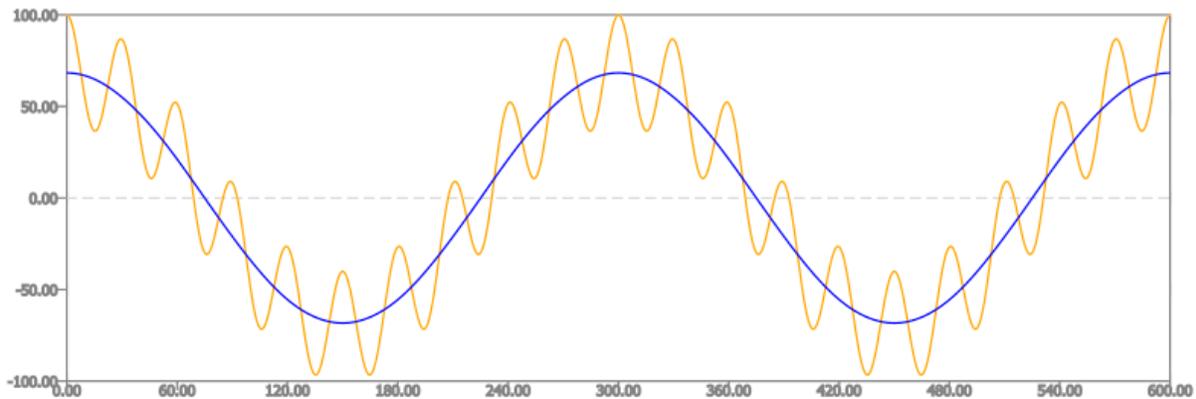

Figure 3.2 Original signals (in orange) and the first residue function (in blue)

## 4   Data filtering

The IMF components decomposed with the derivative initialization method reveal all the riding waves, on each oscillation cycle of the IMF, no further riding wave can be separated. This result can serve to do the data filtering on the IMF components directly by removing the unwanted cycles.

The filtering conditions can be the jump time on the increasing or decreasing edges of the IMF or the amplitude of the extrema on the IMF. The jump time limit relates to the frequency filtering, while the amplitude limit serves to filter out the low power signals. By blocking the short jump time signal, a low pass filter is obtained, by blocking the signal in one range or multiple ranges of jump time, a band block or multiple bands block filter follows. The filtered data function can be reconstructed by combining the passed IMF and the residue function.

## 4.1   Locate the extrema to be filter out and construct the passed function

The filtering start by checking all the extrema of the IMF, if the jump time of the front edge or the back edge of an extremum falls in the given limits, or if the absolute value of the extremum is non zero and less than the amplitude limit, the extremum will be marked for filtering out. By combining the contiguous marked extrema into a single list, all the marked extrema will yield disjoint lists.

Calculate the median points for all the marked IMF extrema using formula 3.2 if the amplitude of the extremum is non-zero, in this case $E(t_i)$ is the value on the IMF, otherwise the amplitude is zero and the median point is identical to the extremum.

For each list of marked extrema, use the median points as the control points and add two extra control points on both ends: the extrema located before the first point in the list and the extrema after the last point in list. A cubic spline function is constructed with this control points with the condition of first derivative equals zero on both ends.

The passed function is glued together with the spline for the marked lists and the original IMF elsewhere.

The process repeats until no extrema is marked or the maximum difference between the successive passed functions is less than a given limit.

## 4.2   Filtering algorithm

The filtering is done through the following steps:

1) Use the method in chapter 3 to find one component of IMF and residue function for the input data to be filtered;
2) Mark the extrema on the IMF to be filtered out by the given jump time or the amplitude limits, if no marked extrema found, the filtering is complete and go to step 5). Otherwise build the marked lists;
3) With the marked lists, calculate the control points and build the cubic spline for each list;
4) Construct the passed function by gluing the spline of each marked list and the IMF elsewhere. Calculate the filtered data function as the sum of the passed function and the residue function. If the maximum point difference to the previous version of the passed function is less than the given limit, the filtering is complete and proceeds to next step. Otherwise use the filtered data function as the input data and go back to step 1);
5) Calculate the blocked function as difference between the original input data and the filtered data.

# 5   Numerical example

The results in the following example are generated with the demo program coded in MS Visual C# 2008. The source code is distributed under BSD open source license.

## 5.1    Original data and last residue

The simulated data is randomly generated. The figure below displays the original data function (in gray) and the last residue function (in blue); the residue reflects the general trend of the data.

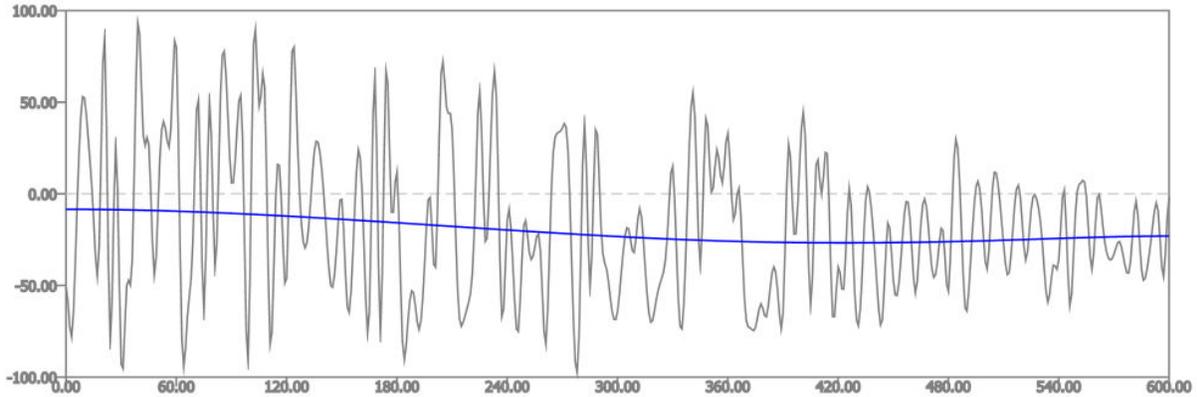

## 5.2    Intrinsic modes function components

The following figures display the six IMF components. The IMF reflects the fluctuation of the data on different time scales.

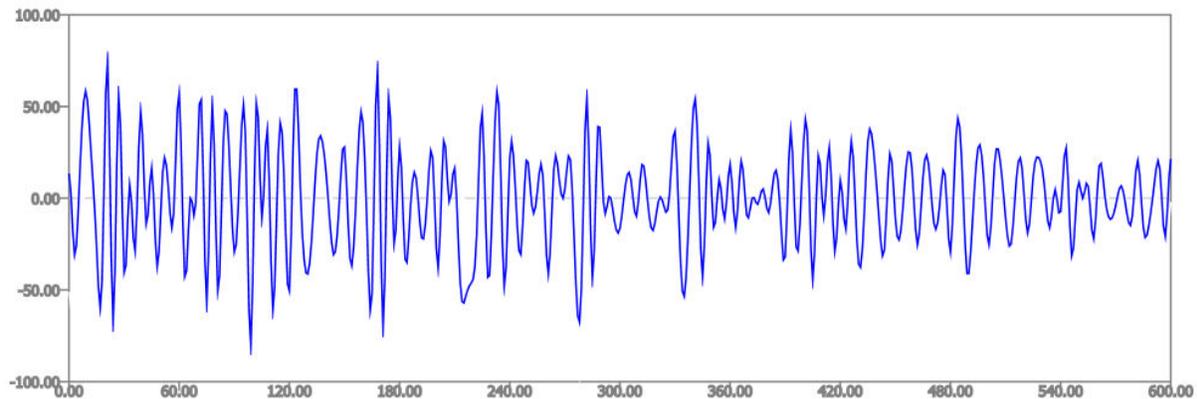

IMF component 1, Extrema count: 147, Value range: [-80.422, 74.992], Iterations: 3, Delta: 2.106113 at 368.00

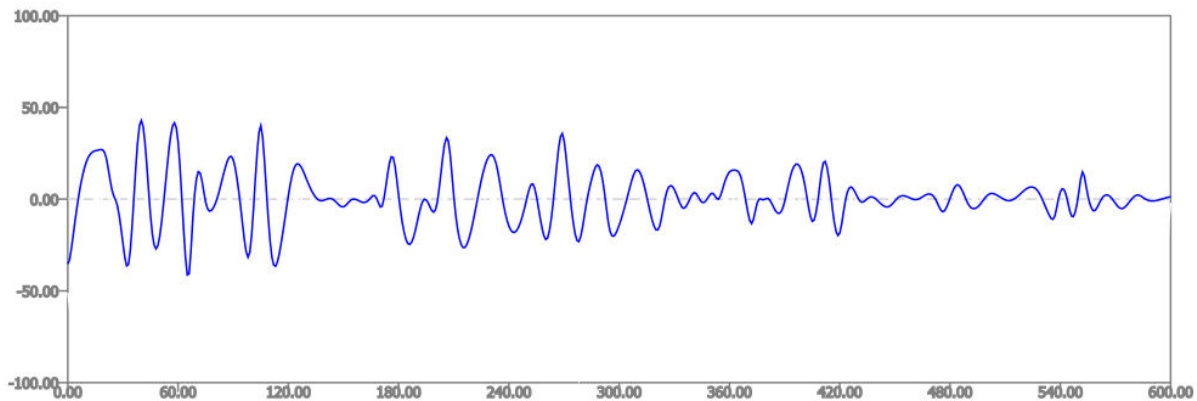

IMF component 2, Extrema count: 76, Value range: [-41.297, 42.876], Iterations: 2, Delta: 10.336595 at 20.00

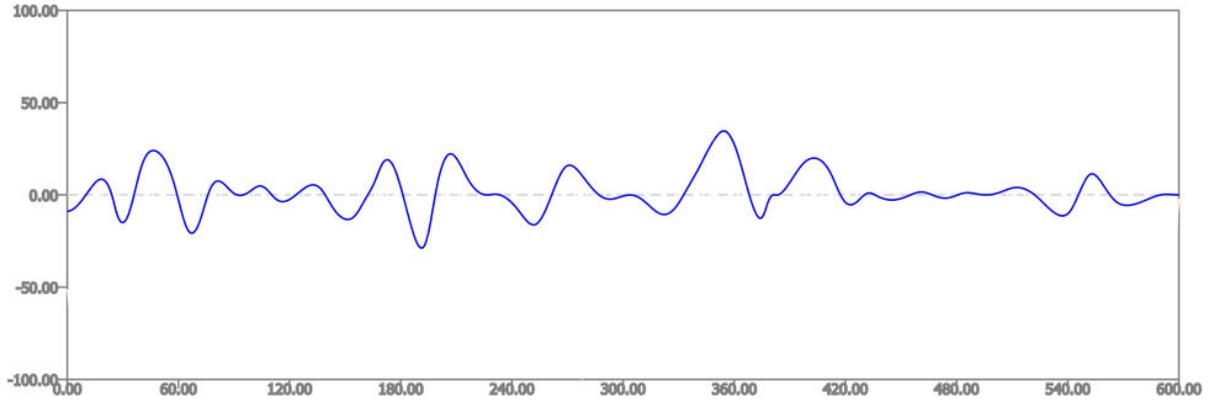

IMF component 3, Extrema count: 37, Value range: [-28.815, 34.694], Iterations: 8, Delta: 0.279790 at 37.00

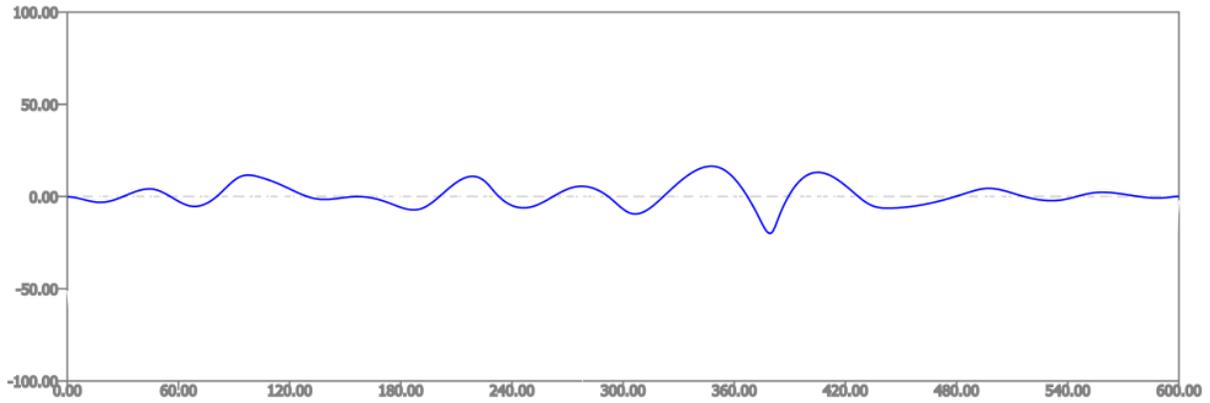

IMF component 4, Extrema count: 19, Value range: [-19.965, 16.451], Iterations: 7, Delta: 0.000000 at 178.00

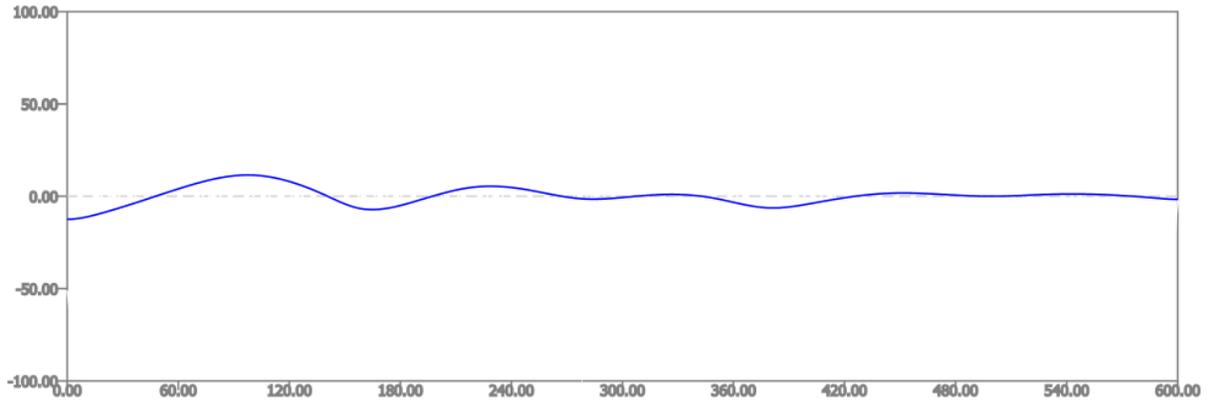

IMF component 5, Extrema count: 9, Value range: [-12.433, 11.527], Iterations: 4, Delta: 0.000000 at 381.00

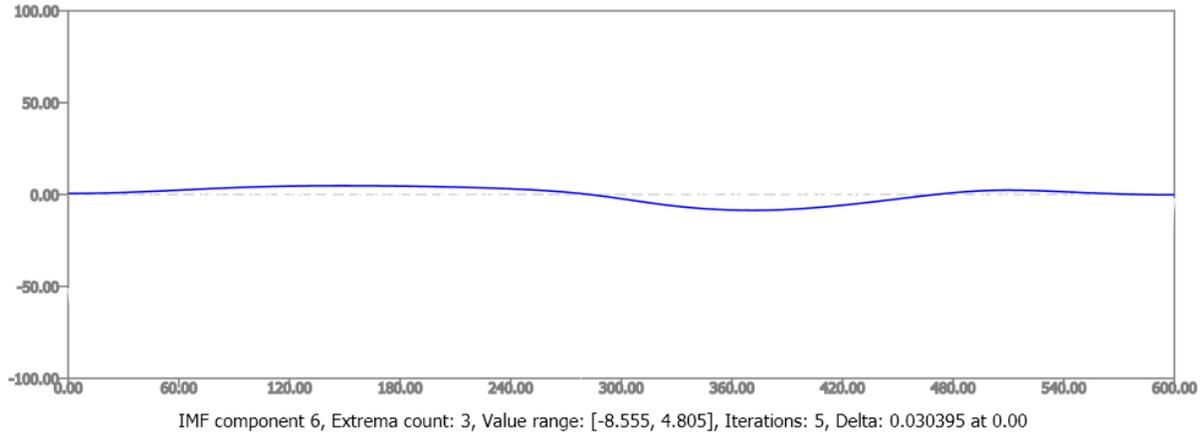

IMF component 6, Extrema count: 3, Value range: [-8.555, 4.805], Iterations: 5, Delta: 0.030395 at 0.00

## 5.3    The residue functions

The generated figures below display the six residue components (in blue), the input data or the residue of the higher frequency mode (in orange).

The residue of different mode reflects the trend of the data changes in different time scales.

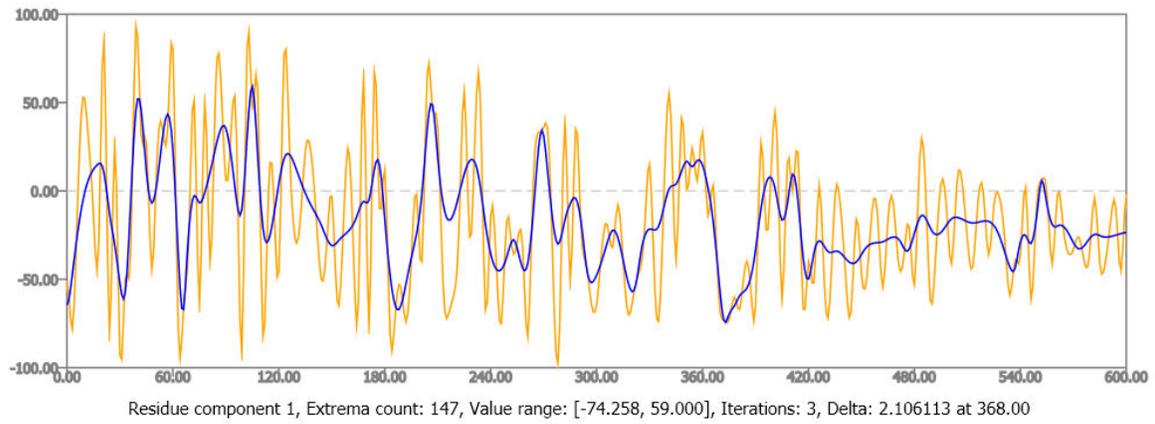

Residue component 1, Extrema count: 147, Value range: [-74.258, 59.000], Iterations: 3, Delta: 2.106113 at 368.00

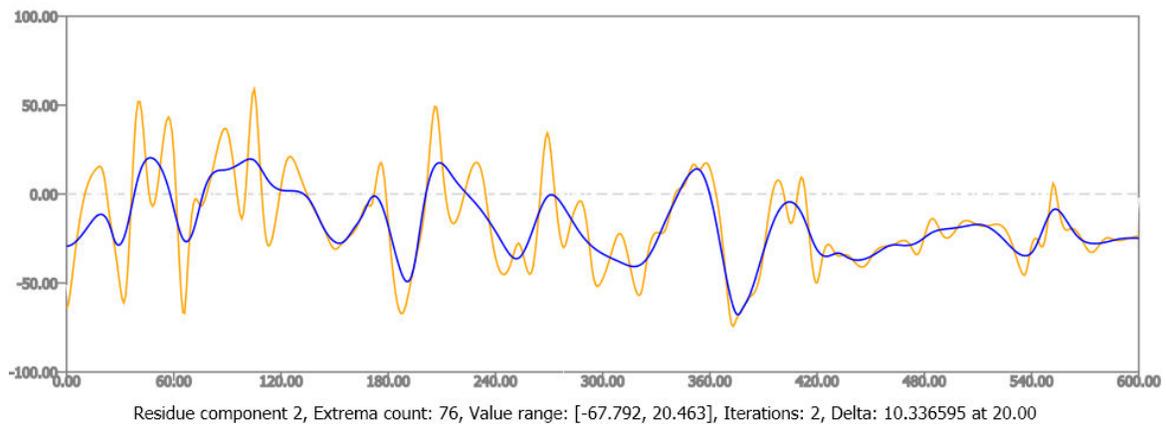

Residue component 2, Extrema count: 76, Value range: [-67.792, 20.463], Iterations: 2, Delta: 10.336595 at 20.00

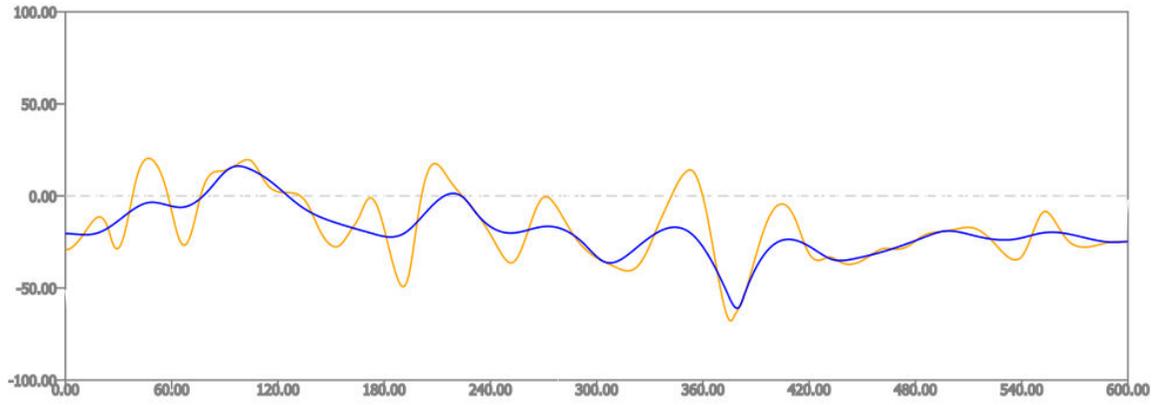

Residue component 3, Extrema count: 37, Value range: [-60.963, 16.232], Iterations: 8, Delta: 0.279790 at 37.00

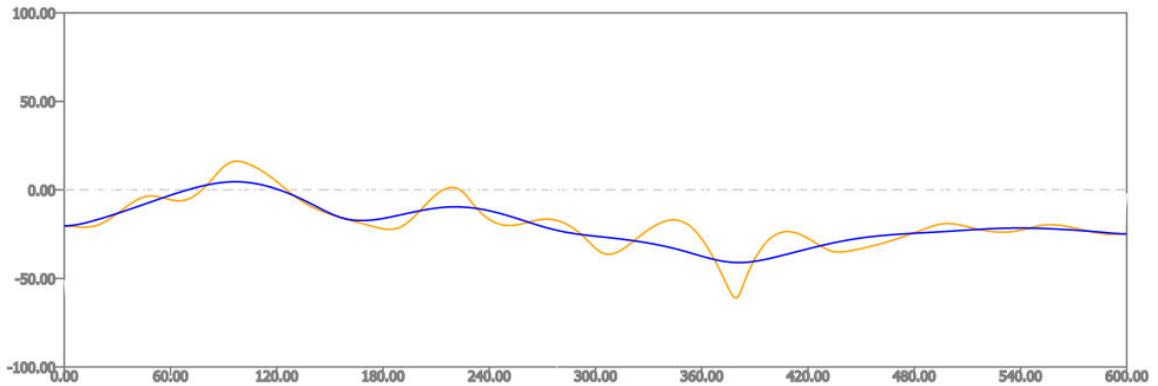

Residue component 4, Extrema count: 19, Value range: [-41.035, 4.610], Iterations: 7, Delta: 0.000000 at 178.00

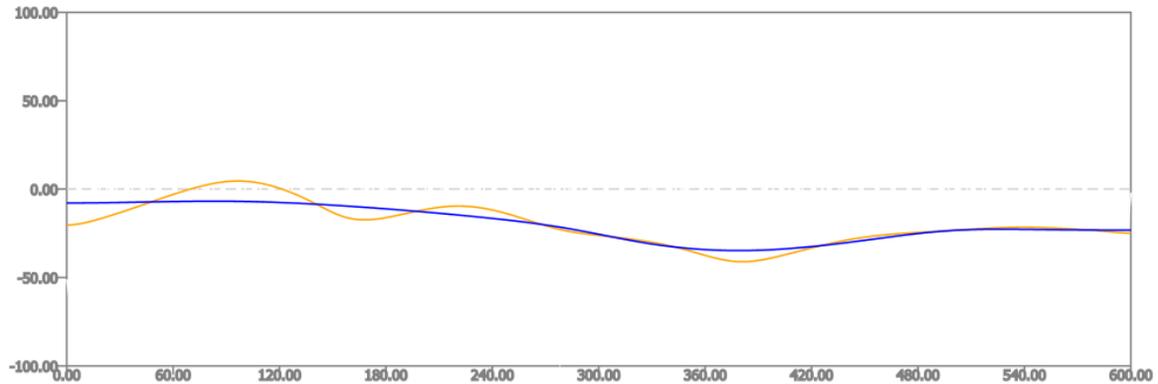

Residue component 5, Extrema count: 9, Value range: [-34.742, -6.827], Iterations: 4, Delta: 0.000000 at 381.00

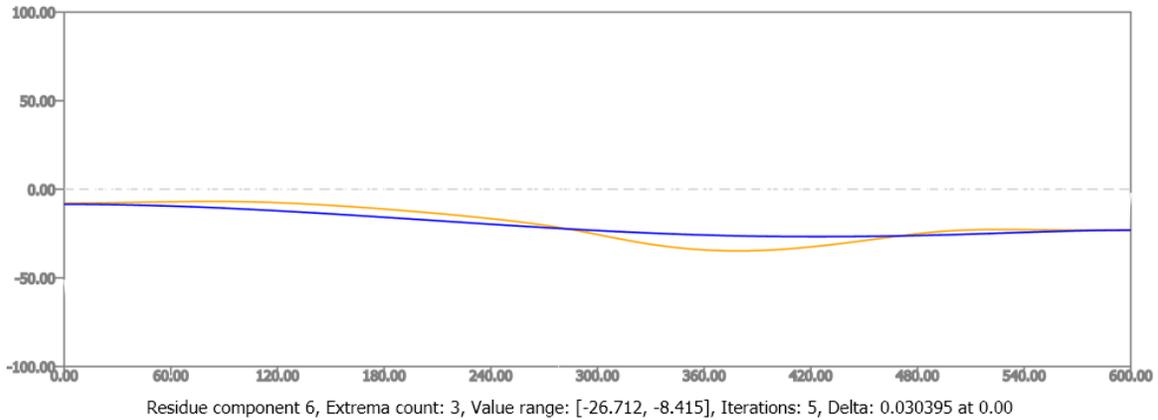

Residue component 6, Extrema count: 3, Value range: [-26.712, -8.415], Iterations: 5, Delta: 0.030395 at 0.00

## 5.4   Comparison of different ways to initialize the IMF

The following two figures display the IMF component with different initialization methods on the same input data, the first one is initialized with the extrema of the first derivative, the second is initialized directly with the extrema of the data function.

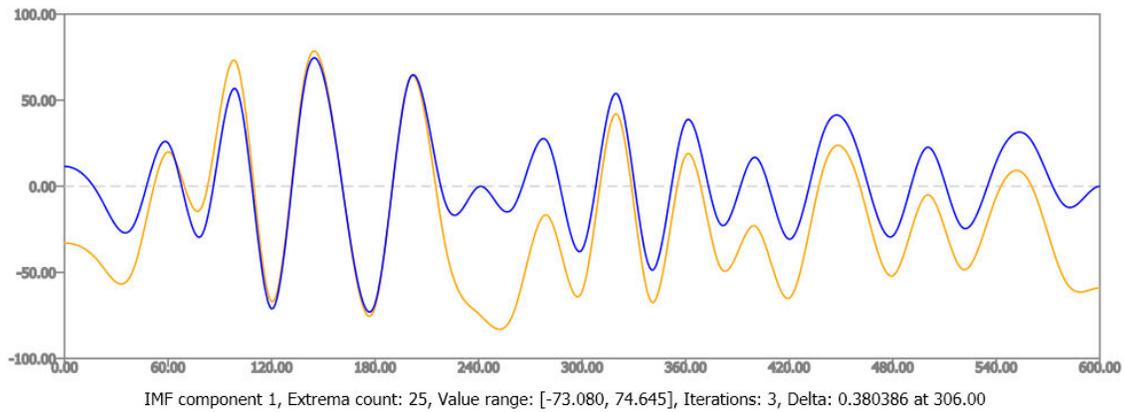

IMF component 1, Extrema count: 25, Value range: [-73.080, 74.645], Iterations: 3, Delta: 0.380386 at 306.00

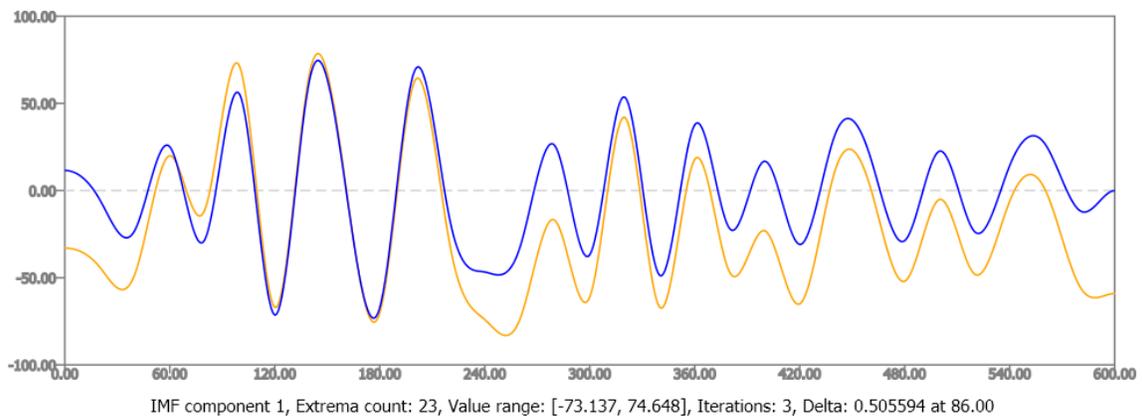

IMF component 1, Extrema count: 23, Value range: [-73.137, 74.648], Iterations: 3, Delta: 0.505594 at 86.00

The riding wave located around 240 is separated in the first figure, but not in the second one. In some application, if the riding wave is undesired noise or ignorable fluctuation, it may be beneficial to use the second IMF initialization method. But the data filtering method is always better approach for noise removal.

5.5   Data filtering example

The figure below displays the passed IMF (in blue) and filtered data (in brown) with the condition of jump time less than 20 for the same input data as in previous section (5.4). The riding waves around 240 and in between 360 to 421 are apparently filtered out.

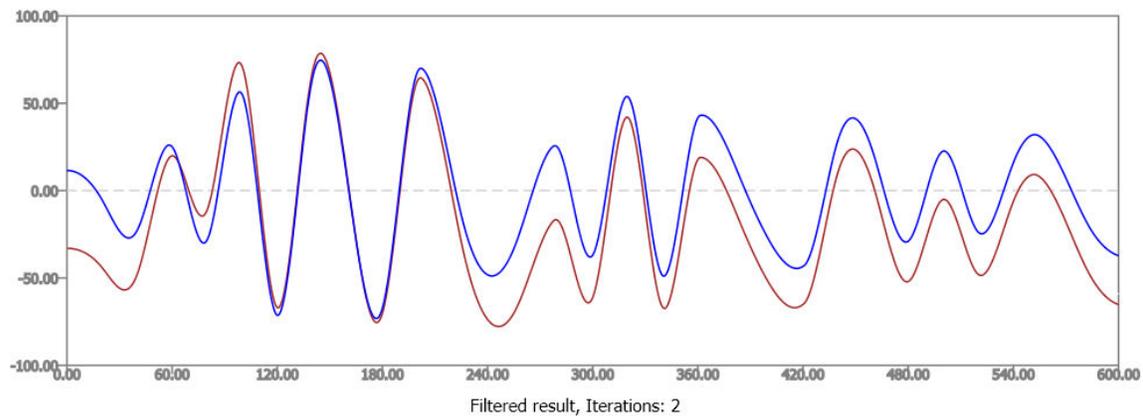

Filtered result, Iterations: 2

# 6   Conclusion

The fast IMD presented in this paper provides a fast and reliable way to decompose the time series data into a unique set of intrinsic mode components. It resolves the problems of convergence and performance related with the EMD method. It also fixes the leaking problem of high frequency ripples into the residue with the sawtooth method. An efficient data filtering method on frequency and amplitude of the IMF is introduced with the similar algorithm, the resulting filter has very little side effect on the useful data. This method will serve as a powerful tool to separate the nonlinear and non-stationary time series data into the trend (residue function) and the oscillation (IMF) on different time scales, and will find application in many fields where traditionally Fourier analysis method or Wavelet analysis method dominate.

**Reference**

[1] The Empirical Mode Decomposition and the Hilbert Spectrum for Nonlinear and Nonstationary Time Series Analysis, Proceedings of the Royal Society of London, A (1998) v. 454, 903-995, Norden E. Huang, Zheng Shen, Steven R. Long, Manli C. Wu, Hsing H. Shih, Quanan Zheng, Nai-Chyuan Yen, Chi Chao Tung and Henry H. Liu

[2] HHT Tutorial, Presentation, July 21, 2004, Norden E. Huang

[3] The Hilbert Transform, Mathias Johansson

[4] Hilbert Huang Transform Technologies, http://www.fuentek.com/technologies/hht.htm

[5] Fast Intrinsic Mode Decomposition of Time Series Data with Sawtooth Transform, Louis Y. Lu, http://arxiv.org/ftp/arxiv/papers/0710/0710.3170.pdf

[6] A Fast Empirical Mode Decomposition Technique for Nonstationary Nonlinear Time Series, Christopher D. Blakely